\documentclass[10pt, a4paper]{article}
\usepackage[english]{babel}
\usepackage{amsmath}
\usepackage{cmll}
\usepackage{amssymb}
\usepackage{amsthm}
\usepackage{bbm}
\usepackage{mathrsfs}
\usepackage{stmaryrd}
\usepackage[all]{xy}
\usepackage{bussproofs}

\newcommand{\sse}{\leftrightarrow}
\newcommand{\arr}{\longrightarrow}

\newcommand{\imply}{\rightarrow}
\newcommand{\C}{\mathbb{C}}

\newcommand{\cl}{\textbf{clop}}

\newcommand{\CH}{\textsf{\textbf{KH}}}

\newdir{|>}{%
!/4.5pt/@{|}*:(1,-.2)@^{>}*:(1,+.2)@_{>}}

\mathcode`\:="603A     
\mathchardef\colon="303A  

\mathcode`\<="4268     
\mathcode`\>="5269     
\mathchardef\gt="313E  
\mathchardef\lt="313C  

\theoremstyle{definition}

\date{}
\author{Fabio Pasquali}

\begin{document}
\title{A tripos based on compact Hausdorff spaces}
\maketitle
In these short notes we show that the category \CH\ of compact Hausdorff spaces is the base of a tripos. There are several consequences of this fact which we shall not discuss here and they come from the well established theory of triposes (one relevant is that \CH\ is mapped into an elementary topos with a finite limits preserving functor). We address interested readers to \cite{Tripos,tripinret}. 
\subsection*{\hfil Triposes\hfil}
Let \textbf{Hyt} be the category of Heyting algebras and Heyting homomorphisms. If $L$ is in \textbf{Hyt} and $x,y$ are in $L$, we shall denote by $x\wedge y$ their meet and by $x\imply y$ the Heyting implication. We use $x\sse y$ as the usual abbreviation for $x\imply y \wedge y\imply x$.\\\\
The following definition of tripos is deduced from \cite{tripinret}.\\\\
A \textit{tripos} is a pair $(\C,P)$ where $\C$ is a category with finite products and $P$ a functor $$P:\C^{op}\arr \textbf{Hyt}$$
such that
\begin{itemize}
\item[i)] for each arrow $f: A\arr B$ the homomorphism $$P(f):P(B)\arr P(A)$$ has a left adjoint $\exists_{f}$ and a right adjoint $\forall_{f}$ satisfying the Beck-Chevalley condition, i.e. for every pullback of the form 
\[\xymatrix{
X\times\Gamma\ar[d]_-{id_X\times k}\ar[r]^-{\pi_\Gamma}&\Gamma\ar[d]^{k}\\
X\times\Delta\ar[r]_-{\pi_\Delta}&\Delta
}\]
it is $\exists_{\pi_\Gamma} (id_X\times k)^*= k^*\exists_{\pi_\Delta}$ and $\forall_{\pi_\Gamma} (id_X\times k)^*= k^*\forall_{\pi_\Delta}$.
\item[ii)] $(\C,P)$ has \textit{weak power objects}, i.e. for every $X$ in $\C$, there exists $\mathbb{P}X$ in $\C$ and a formula $\in_X$ in $P(X\times \mathbb{P}X)$ such that for every object $Y$ in $\C$ and every formula $\gamma$ in $P(X\times Y)$ there exists $\{\gamma\}:Y\arr \mathbb{P}X$ such that $(id_X\times \{\gamma\})^*\in_X = \gamma$.
\end{itemize}
There are many redundancies in the previous definition \cite{Tripos}. We only mention the following.\\\\
For every $X$ in $\C$ it is $$\delta_X=\forall_{<\pi_1,\pi_2>}(<\pi_1,\pi_3>^*\in_X\sse<\pi_2,\pi_3>^*\in_X)$$
where $\pi_i$, whih $i=1,2,3$ are projections from $X\times X\times \mathbb{P}X$ to each of the factors.
\subsection*{\hfil A \CH-based tripos\hfil}
We denote by \CH, the full subcategory of \textbf{Top} on Compact Hausdorff spaces. All the facts on general topology which are not proved are from \cite{Kelley}\\\\
Consider the pair $(\CH, \cl)$ where $\cl:\CH^{op}\arr\textbf{Hyt}$
is the functor that maps every Compact Hausdorff space into the Boolean algebra of its clopen sets and every continuous function into the inverse image functor.\\\\
The binary product of compact Hausdorff spaces is compact Hausdorff. The projections are always open maps. If the two factors of the product are compact, the projections are also closed. Then projections in \CH\ are clopen maps. Therefore, if $\pi_Y:X\times Y\arr Y$ is a projection in \CH, then $$\text{Im}_{\pi_Y}:\cl(X\times Y)\arr\cl(Y)$$ is left adjoint to $\pi_Y^{-1}$. Since $\cl(Y)$ is a Boolean algebra the map $\pi_Y^{-1}$ has also a right adjoint $\forall_{\pi_Y} = \neg \text{Im}_{\pi_Y}\neg$ where $\neg$ is the operation of taking complements.\\\\
It remains to prove that $(\CH, \cl)$ has weak power objects.\\\\
For every $Y$ in \CH\ and every $A$ in $\cl(Y)$, $A$ equipped with the subspace topology is in $\CH$, since subspaces of Hausodrff spaces are Hausdorff and closed subspaces of compact spaces are compact.\\\\
The following square is a pullback in $\CH$
\[
\xymatrix{
A\ar[d]_-{\lfloor A\rfloor}\ar[r]&\textbf{1}\ar[d]^{t}\\
Y\ar[r]_-{\chi_A}&\textbf{2}
}
\]
where \textbf{1} and \textbf{2} are $\{1\}$ and $\{0,1\}$ endowed with the discrete topology, $t:\textbf{1}\arr\textbf{2}$ is the constant map to $1$ and $\chi_A$ is the characteristic function of $A$, which is continuous since $A$ is clopen.
\\\\
For every topological space $B$ we denote its Alexandroff compactification by $B_\infty$. The set of points $|B_\infty |$ is $B\cup\{\infty\}$ where $\infty$ does not belongs to $|B|$. The open sets of $B_\infty$ are the open sets of $B$ together with all subsets $V$ of $|B_\infty|$ such that $\infty$ is an element of $V$ and $B_\infty\setminus V$ is closed and compact in $B$.\\\\The inclusion $i:B\arr B_\infty$ is an open map. Moreover $B_\infty$ has the property that if $f:A\times B \arr \textbf{2}$ is continuous and $B$ is in \CH, then the extension $f_\infty: A\times B_\infty \arr \textbf{2}$ determined by
\[
f_\infty(a,b) =
  \begin{cases}
   f(a,b) & \text{if } (a,b)\ \text{is in}\ A\times B\\
   0       & \text{otherwise }
  \end{cases}
\]
is continuous. In fact $f_\infty^{-1}\{1\} = i[f^{-1}\{1\}]$ is open since $f^{-1}\{1\}$ is open and $i$ is an open inclusion. Analogously $i[f^{-1}\{0\}]$ is open. On the other and $$f_\infty^{-1}\{0\}= i[f^{-1}\{0\}]\cup A\times\{\infty\}$$
The singleton $\{\infty\}$ is open in $B_\infty$ as $B_\infty\setminus\{\infty\}=B$ which is closed and compact.\\\\
For every $A$ and $B$ in \CH\ and for every $\phi$ in $\cl(A\times B)$ consider its characteristic function $$\chi_\phi:A\times B\arr \textbf{2}$$
Since \textbf{2} and $A$ are compact and Hausdorff the space $\textbf{2}^A$ endowed with the compact-open topology is an exponential in \textbf{Top} (but not necessary in $\CH$). Then there are continuous functions $$ev: A\times \textbf{2}^A\arr\textbf{2}$$
$$\overline{\chi_\phi}: B\arr \textbf{2}^A$$ such that $ev(id_A\times \overline{\chi_\phi}) = \chi_\phi$. Now consider the diagram
$$
\xymatrix{
A\times (\textbf{2}^A)_\infty\ar[rrd]^-{ev_\infty}&&\\
A\times \textbf{2}^A\ar[rr]^-{ev}\ar[u]^-{id_A\times i}&&\textbf{2}\\
A\times B\ar[rru]_-{\chi_\phi}\ar[u]^-{id_A\times\overline{\chi_\phi}}&&
}
$$
We claim that $(\textbf{2}^A)_\infty$ is a weak power object of $A$. Define $$\in_A = ev_{\infty}^{-1}(\{1\})$$ For every clopen set $\phi$ in $A\times B$ define $\{\phi\} = i\overline{\chi_{\phi}}$. Then $$(id_A \times \{\phi\})^{-1}\in_A =(id_A\times\overline{\chi_\phi})^{-1}(id_A\times i)^{-1}ev_\infty^{-1}\{1\}= \chi_\phi^{-1}\{1\}=\phi$$
Therefore the space $(\textbf{2}^A)_\infty$ is a weak power object of $A$ provided that it is an object of \CH. The space $\textbf{2}^A$ is Hausdorff since \textbf{2} is Hausdorff. For every $f$ in $\textbf{2}^A$ the sets $f^{-1}\{1\}$ and $f^{-1}\{0\}$ are compact as they are closed subsets of a compact space. Then the sets of all continuous functions $\mathcal{C}(f^{-1}\{1\},\{1\})$ and $\mathcal{C}(f^{-1}\{0\},\{0\})$ are open in the compact-open topology, therefore $$\{f\}=\mathcal{C}(f^{-1}\{1\},\{1\})\cap\mathcal{C}(f^{-1}\{0\},\{0\})$$
is open, which implies that $\textbf{2}^A$ is discrete. Discrete Hausdorff spaces are locally compact. The Alexandroff compactification of a space is Hausdorff exactly when the space is locally compact Hausdorff \cite{Kelley}, from which the claim.

\end{document}